\documentclass[11pt]{article}
\usepackage{graphicx,setspace}
\usepackage{amsmath,amssymb,amsthm}
\newcommand{\be}{\begin{equation}}
\newcommand{\ee}{\end{equation}}
\newcommand{\bea}{\begin{eqnarray}}
\newcommand{\beas}{\begin{eqnarray*}}
\newcommand{\no}{\nonumber}
\newcommand{\eea}{\end{eqnarray}}
\newcommand{\eeas}{\end{eqnarray*}}

\newcommand{\prob}[1]{\mathsf{Pr}\left( #1 \right)}

\newcommand{\EXP}[1]{\mathsf{E}\!\left(#1\right) }
\newcommand{\COV}[1]{\mathsf{Cov}\left( #1 \right)}
\newcommand{\VAR}[1]{\mathsf{VAR}\!\left(#1\right) }
\newcommand{\remove}[1]{}

\newtheorem{theorem}{\bf Theorem}

\newcounter{cnt1}

\newcounter{cnt3}
\newcommand{\blr}{\begin{list}{$($\roman{cnt1}$)$}
 {\usecounter{cnt1} \setlength{\topsep}{0pt}
 \setlength{\itemsep}{0pt}}}
\newcommand{\bla}{\begin{list}{$($\betaph{cnt2}$)$}
 {\usecounter{cnt2} \setlength{\topsep}{0pt}
 \setlength{\itemsep}{0pt}}}
\newcommand{\bln}{\begin{list}{$($\arabic{cnt3}$)$}
 {\usecounter{cnt3} \setlength{\topsep}{0pt}
 \setlength{\itemsep}{0pt}}}
\newcommand{\el}{\end{list}}

\def\half{\frac{1}{2}}

\def\rar{\rightarrow}

\newcommand{\ep}{\epsilon}

\newcommand{\del}{\delta}
\newcommand{\lam}{\lambda}

\newcommand{\sg}{\sigma}

\def\cC{{\mathcal C}}
\def\us{{\underline s}}
\begin{document}
\begin{titlepage}
\begin{center}
{\bf Limit Laws for $k$-Coverage of Paths by a Markov-Poisson-Boolean Model} \\
\vspace{0.2in} {Srikanth K. Iyer \footnote{corresponding author: skiyer@math.iisc.ernet.in}$^,$\footnote{Research Supported in part by UGC SAP -IV and DRDO grant No. DRDO/PAM/SKI/593}}\\
Department of Mathematics,
Indian Institute of Science, Bangalore, India. \\
\vspace{0.1in}
D. Manjunath \footnote{Work carried out in the Bharti Centre for
       Communications Research at IIT Bombay and supported in
       part by a grant from the Ministry of Information Technology,
       Government of India.}\\ Department of Electrical Engineering,
Indian Institute of Technology Bombay, Mumbai, India. \\
\vspace{0.1in}
D. Yogeshwaran \footnote{Supported in part by a grant from EADS, France.}\\
INRIA/ENS TREC, Ecole Normale Superieure, Paris, France.
\end{center}
\vspace{0.2in}
%
\sloppy
\begin{center} {\bf Abstract} \end{center}

\begin{center} \parbox{4.8in}
{Let $P := \{ X_i\}_{i \geq 1}$ be a stationary Poisson point process
  in $\Re^d,$ $\{C_i\}_{i \geq 1}$ be a sequence of i.i.d. random sets
  in $\Re^d,$ and $\{Y_i^t; \; t \geq 0 \}_{i \geq 1}$ be i.i.d.
  $\{0,1\}$-valued continuous time stationary Markov chains.  We
  define the \emph{Markov-Poisson-Boolean model} $\cC_t := \{Y_i^t(X_i
  + C_i), i \geq 1 \}.$ $\cC_t$ represents the coverage process at
  time $t.$
  We first obtain limit laws for $k$-coverage of an area at an
  arbitrary instant. We then obtain the limit laws for the
  $k$-coverage seen by a particle as it moves along a one-dimensional
  path.  } \\
\vspace{0.2in} \today
\end{center}

\vspace{0.2in}
{\sl AMS 1991 subject classifications}: \\
\hspace*{0.5in} Primary:   60D05, 60G70\\
\hspace*{0.5in} Secondary:  05C05, 90C27\\
{\sl Keywords:} Poisson-Boolean model, coverage, Markov process,
sensor networks, target tracking.

\end{titlepage}

\remove{

\title{Limit Laws for $k$-Coverage of Paths by a Markov-Poisson-Boolean Model}

\author{ Srikanth.~K.~Iyer$^{\dagger}$, D.~Manjunath$^{\ddag}$ and
  D.~Yogeshwaran$^{\sharp}$}
\thanks{ $\dagger$ Contact Author:
    Departement of Mathematics, Indian Institute of Science Bangalore
    - 560012, INDIA. \texttt{skiyer@math.iisc.ernet.in}} \thanks{
    $\ddag$ Dept. of Electrical Engineering, Indian Institute of
    Technology Bombay, Mumbai - 400076, INDIA.
    \texttt{dmanju@ee.iitb.ac.in}} \thanks{ $\sharp$ INRIA/ENS TREC,
    Ecole Normale Superieure, Paris-75005 , FRANCE.
    \texttt{yogesh@di.ens.fr} }

\date{\today}
\maketitle
\begin{abstract}

  Let $P := \{ X_i\}_{i \geq 1}$ be a stationary Poisson point process
  in $\Re^d,$ $\{C_i\}_{i \geq 1}$ be a sequence of i.i.d random sets
  in $\Re^d,$ and $\{Y_i^t; \; t \geq 0 \}_{i \geq 1}$ be i.i.d.
  $\{0,1\}$-valued continuous time stationary Markov chains.  We
  define the \emph{Markov-Poisson-Boolean model} $\cC_t := \{Y_i^t(X_i
  + C_i), i \geq 1 \}.$ $\cC_t$ represents the coverage process at
  time $t.$
  We first obtain limit laws for $k$-coverage of an area at an
  arbitrary instant. We then obtain the limit laws for the
  $k$-coverage seen by a particle as it moves along a one-dimensional
  path.

\end{abstract}

\vspace{0.5in}

\noindent Keywords: Poisson-Boolean model, coverage, Markov process,
sensor networks, target tracking.

\newpage
}

\section{Introduction}
\label{sec:Introduction}
\subsection{Motivation}
\label{sec:background}
This paper is motivated by the need to characterize the ability of a
sensor network with randomly deployed unreliable sensor nodes to track
the trajectory of a moving target in the sensor field.  Since the
nodes are randomly deployed, a suitable point process in a
$d$-dimensional metric space, typically $\Re^d$ or $\mathcal{Z}^d,$
can be used to describe the location of the sensor nodes. In this
paper we will assume a homogenous Poisson point process for the distribution of
sensor locations.  A sensor node can detect events or perform
measurements over a `sensing area' or a `footprint.'  The coverage
area of each sensor is described by a suitable sequence of random
$d$-dimensional sets.  Thus analyzing the coverage by the sensor
network involves the analysis of an equivalent \emph{stochastic
  coverage process.} Coverage analysis usually takes the form of
obtaining statistics for the fraction of the volume of a
$d$-dimensional set that is covered by one or more sensors, e.g., as
in \cite{Hall88}.

Tracking a moving target by a sensor network involves trajectory
estimation from a sequence of position estimates of the target.  The
quality of the trajectory estimate will depend on the parts of the
trajectory that are `covered' by $k$ or more sensors with the value
of $k$ being determined by the estimator used.
%

The sensor nodes are assumed to be unreliable in the sense that they toggle
between two states---`available' and `not available.'  Hence some of
the sensors that are covering the target as it is moving along the
path could become unavailable during the coverage of the target. There
are many reasons for this.  A sensor could be turned off for energy
saving or even energy restoration in the battery.  Alternatively, a
sensor may have temporarily failed. It could also be that the radio
environment is such as to prevent a sensor from communicating with its
neighbors, effectively making it unavailable for the sensing process.
Since the nodes are unreliable, the coverage of the path during the
motion of the target is a random process in time that is determined by
the switching of the sensor states that could have covered the path.
Thus we need to explicitly model the temporal behavior of the
induced-coverage process of the path. We model this ``on-off'' behavior
of the sensors by a two state Markov chain.

The aim of this paper is two-fold. We first extend the asymptotic coverage
results in  \cite{Hall88} to the case of general $k$-coverage.
Asymptotic properties of the covered fraction are derived as the intensity of the
Poisson point process becomes large and the coverage area of the individual
sensors scaled down in such a way that the limiting process covers a non-trivial
 fraction of the operational area. The second objective is to
show how the dynamics of the on-off process affects the coverage
of a linearly moving target.
%
%
\subsection{The Markov-Poisson-Boolean Model}
\label{sec:Markov-Boolean}
Let $P := \{ X_i, i =1,2,\ldots \}$ be a
stationary Poisson process in $\Re^d$ of intensity $\lambda > 0$.  The points of
$P$ can be thought of as locations of sensors in a random sensor
network. Let $C$ be a random closed set in $\Re^d$ independent of $P$ and
having an arbitrary distribution. Throughout the paper we assume that for some fixed $\tau > 0,$ $C
\subset B_{0,\tau},$ almost surely, where $B_{0,\tau}$ is a closed ball
of radius $\tau$ centered at the origin. We will also assume that $\|C\| >0$ with probability 1, where $\| \cdot \|$ denotes
the Lebesgue measure. Let
$\beta = \EXP{\|C\|}.$ Let $C_i$ be i.i.d. copies of $C$.
As in \cite{Hall88}, the $C_i$'s will be called shapes to
distinguish them from the sets $(X_i + C_i)$'s that denote the areas in
$\Re^d$ that are covered by the sensors.
The coverage process $\mathcal{C} \equiv \{ X_i + C_i, i \geq 1 \} $
is called the \emph{Poisson-Boolean model} \cite{Hall88}. Now let
$Y= \{ Y^t, \; t \geq 0\}$ be a $\{0,1\}$-valued continuous time stationary
Markov process independent of all other random variables. Let
$\{Y^t_i, t \geq 0 \}, i \geq 1$ be i.i.d.  copies of $Y.$ $Y_i^t=1$ can be
interpreted to mean that sensor $i$ is `on' at time $t$ and is
available for sensing and $Y_i^t=0$ means that it is `off' at time $t$
and not available for sensing. Define the \emph{Markov-Poisson-Boolean
  model} $\cC_t := \{Y_i^t(X_i + C_i), i \geq 1 \}.$ $\cC_t$
represents the coverage process by the available sensors at time
$t.$ Let $R \subset \Re^d$ be an arbitrary Borel set.
$R$ could correspond to the operational area of the sensor
network. Now, consider a straight~line path $\mathcal{L}$ in $R$ of
length $L$ units, and an object moving along it with a constant
velocity of $c$ units per second. The object starts moving at time
$t=0.$ For any positive integer $k$, let $\zeta_k(t)$ be the
indicator function for the object being `covered,' to be made more
precise later, by $k$ or more sensors at time $t.$ The objective of
this paper is to characterize the random variables $\zeta_k(t)$ and
$\int_0^T\zeta_k(t)dt$, where $T= L/c.$
In the rest of the paper, we refer to the \emph{Poisson-Boolean model}
and the \emph{Markov-Poisson-Boolean model} as the PB and the MPB
models respectively.
\subsection{Previous Work}
\label{sec:prev-work}
Coverage with reliable nodes, where nodes are always in the
`available' state, has been elaborated in both applied mathematics and
sensor network literature. See \cite{Hall88} for a good comprehensive
first study and \cite{Liu04} for some recent results for the case
$k=1.$ In \cite{Zhang04}, $R$ is a square of side-length $\ell$ and $C$
is a circle of unit area.  It is shown for this PB model, that
if $\lambda$ is given by
\begin{equation}
  \label{eq:k-coverage-zhang}
  \lambda = \log \ell^2 + (k+1) \log \log \ell^2 + c_{\ell},
\end{equation}
and if $c_{\ell} \to \infty$ as $\ell \to \infty$, then $R$ is almost surely
asymptotically $k$-covered.

Asymptotic coverage by unreliable sensor networks at an arbitrary epoch
has been studied in \cite{Shakkottai03,Kumar04}.  In the analysis of
coverage by an unreliable sensor network at an arbitrary epoch, the
stationary probability of a sensor node not being available
essentially `thins' the original deployment process and a standard
analysis with the thinned process applies.  However, for applications
like target tracking or intruder detection, we need to know the
behavior of the coverage process \emph{during} the movement of the
target. When the sensor nodes are unreliable, a node that was sensing
the object may switch from being available to becoming unavailable, or
vice versa. This implies that the coverage of a given point is not
independent either in space or in time.  Thus we need to consider the
dynamics of the transitions from availability to non-availability of
the sensors in the spatio-temporal analysis of the coverage of the
path. The coverage of a line by a two-dimensional PB model was
investigated in \cite{Manohar06, Sundhar06a}.
\subsection{Organization of the Paper and Summary of Results}
\label{sec:orgn-and-summary}
%
In Section~\ref{sec:direct-k-coverage}, we
characterize $k$-coverage in $d$-dimensions. Although
our eventual interest is the characterization of the coverage of a
moving point on a path by the MPB model defined earlier, it is
instructive to first consider the $k$-coverage of $R \subset \Re^d,$
by the MPB model $\cC_t$ at an arbitrary
instant. At an arbitrary
instant, let $V_k(R)$ be the volume of an arbitrary $d$-dimensional
set $R \subset \Re^d$ that is not covered by $k$ or more sensors. We
obtain a strong law and central limit theorem for $V_k(R)$. The proof techniques
are in general similar to those in Chapter 3 of \cite{Hall88}.
In the second part of Subsection~\ref{sec:limit_laws_k_coverage}, we consider
the special case when the coverage areas are discs of fixed radius $r,$
and obtain a strong law of large numbers for the critical radius
required for complete $k$-coverage in dimension $d=2$.

In Section~\ref{sec:markov-boolean}, we consider the MPB model defined
earlier, which is the PB model but now with unreliable sensors. We analyze the path coverage for a linearly
moving target that is in the sensor field for $T$ units of time.
Without loss of generality, let this interval be $(0,T).$ Let
$V_{k,T}$ be the total time in $(0,T)$ that the target is not tracked
by $k$ or more sensors.  For pedagogic convenience we consider $k=1$
and obtain a strong law and central limit theorem for $V_{1,T}.$
The proof techniques of Section \ref{sec:direct-k-coverage} can be
combined with those of Section~\ref{sec:markov-boolean} to extend the
results to the case of $k > 1.$ We have separated $k$-coverage from
the Markovian on-off dynamics to maintain clarity of exposition. Each
of these two components operates independently in the computations and
the expressions involving $k$-coverage are at times lengthy.
%
%
\section{$k$-Coverage}
\label{sec:direct-k-coverage}
\subsection{Preliminaries}
\label{sec:notation-prelims}
For a point $x \in \Re^d$, let $\chi_m(x)$ be the indicator variable
that $x \in X_i + C_i$ for exactly $m$ points in $P,$ i.e.,
\begin{equation}
  \chi_m(x) =
  \begin{cases}
    1 & \mbox{if $x \in X_i + C_i$ for exactly $m$ of $i=1,2,\ldots$} \\
    0 & \mbox{otherwise.}
  \end{cases} \no
\end{equation}
%
%
\begin{eqnarray}
  \EXP{\chi_m(x)} &=& \prob{x \in (X_i+C_i)\mbox{ for exactly $m$
      points}}\nonumber \\
  &=& \prob{\mbox{for exactly $m$ points } X_i \in (x
    - C_i)} \nonumber \\
  &=& \prob{\mbox{for exactly $m$ points }  X_i \in C_i} .
  \label{eq:chi_m-of-x}
\end{eqnarray}
The last equality follows from the stationarity of $P.$ Recall that
$C_i \subset B_{0,\tau}$ and $\beta = \EXP{\| C \|}.$  Given that $N$ points of $P$ lie within
$B_{0,\tau},$ the probability that exactly $m$ of these points cover
the origin is given by
\begin{displaymath}
  \binom{N}{m} \left( 1 - \frac{\beta}{\|B_{0,\tau}\|} \right)^{N-m}
  \left( \frac{\beta}{\|B_{0,\tau}\|} \right)^m .
\end{displaymath}
Since $N$ is Poisson with mean $\lam \|B_{0,\tau}\|$, we obtain
\begin{eqnarray}
  \EXP{\chi_m(x)}   & = & \EXP{ \binom{N}{m} \left( 1 -
      \frac{\beta}{\| B_{0,\tau} \|}
    \right)^{N-m} \left( \frac{\beta}{\| B_{0,\tau} \|} \right)^m }
\;\; = \;\; \frac{e^{-\beta\lam}(\lam\beta )^m}{m!}
  \label{eq:prob-x-covered-by-m}.
\end{eqnarray}
Let $R \in \Re^d$ be a $d$-dimensional set. For $k > 0,$ we define the
$k$-vacancy within $R$, $V_k(R),$ to be the $d$-dimensional volume of
the part covered by at most $k-1$ random sets of $\mathcal{C},$ i.e.,
\begin{equation}
  \label{eq:V_k-defn}
  V_k(R) \equiv \sum_{m = 0}^{k-1} \int_R \chi_m(x) \ dx.
\end{equation}
The indicator variable for the $k$-vacancy of a point $x$ will be
denoted by $V_k(x),$ i.e., $V_k(x) := \sum_{m = 0}^{k-1} \chi_m(x).$
$\| R \| - V_k(R)$ will be called the $k$-\emph{coverage} of $R$.
Since $R$ is fixed throughout the paper, we will omit the reference
to it in the notation and write $V_k(R)$ as $V_k$.

Some of the early derivations mimic that in \cite{Hall88} and we give
it here for the sake of completeness.
From (\ref{eq:prob-x-covered-by-m}) and Fubini's theorem,
\begin{equation}
  \EXP{V_k} \ =\ \sum_{m = 0}^{k-1} \int_R \EXP{\chi_m(x)}dx  =
  \sum_{m = 0}^{k-1} \|R\|\frac{e^{-\lambda \beta}(\lambda
    \beta)^m}{m!}.
  \label{eq:expected-V_k}
\end{equation}
We now derive the variance of $V_k.$ If a point $x$ is covered by
$X_i,$ then $x \in X_i + C_i$ or $X_i \in x - C_i$. Similarly, if $x$
is not covered by $X_i$ then $X_i \in x - C_i^c$ where $C^c = \Re^d \setminus C.$ We
use this to first obtain the probability that two points $x_1$ and
$x_2$ are covered by exactly $m$ and $n$ sensors respectively which is
then used to obtain the variance of $V_k.$ We make the following
observations regarding the location of $X_i$ relative to points $x_1$
and $x_2.$
\begin{itemize}
\item If $X_i$ covers $x_1$ and $x_2$, then $X_i \in B_1^i(x_1,x_2) :=
  (x_1 - C_i) \cap (x_2 - C_i)$. Further, $\|B_1^i\|$ has the same
  distribution as $\|B_1\|$ where $B_1(x_1-x_2) = (x_1 - x_2 + C) \cap
  C.$
\item If $X_i$ covers $x_1$ and not $x_2$, then $X_i \in
  B_2^i(x_1,x_2) := (x_1 - C_i) \cap (x_2 - C^{c}_{i})$ and
  $\|B_2^i\|$ has the same distribution as $\|B_2\|$ where
  $B_2(x_1-x_2) = (x_1 - x_2 + C) \cap C^c.$
\item Similarly, if $X_i$ does not cover $x_1$ but covers $x_2$, then
  $X_i \in B_3^i(x_1,x_2) := (x_1 - C_i^{c}) \cap (x_2 - C_i)$ and
  $\|B_3^i\|$ is equal to $\|B_3\|$ in distribution where
  $B_3(x_1-x_2) = (x_1 - x_2 + C^c) \cap C.$
\end{itemize}
We will suppress the argument of $B_j$ and $B_j^i$ unless required.
Observe that the $B_j$ defined above are mutually disjoint sets.
Further, $\|B_2\|$ and $\|B_3\|$ will have the same distribution.
\begin{eqnarray*}
  \EXP{\chi_m(x_1)\chi_n(x_2)} &=& \prob{\mbox{$x_1$ covered by $m$ \&
      $x_2$ covered by  $n$ sensors}} \\
  &=& \prob{\mbox{$x_1 \in X_i +C_i$ for $m$ sensors and
      $x_2 \in X_i +C_i$ for $n$ sensors}} \\
  &=& \prob{\mbox{$X_i \in x_1 - C_i$ for $m$ sensors and
      $X_i \in x_2 - C_i$ for $n$ sensors}} \\
  & = & \sum_{l=0}^{m \wedge n} \prob{\mbox{$l$ cover $x_1$ \& $x_2$;
      $(m-l)$ cover only $x_1$; $(n-l)$ cover only $x_2$}}.
\end{eqnarray*}
We can proceed as in the derivation of (\ref{eq:chi_m-of-x}) and
(\ref{eq:prob-x-covered-by-m}) and consider a bounded set $A^\prime$ that
contains $(x_1+B_{0,\tau})$ and $(x_2+B_{0,\tau})$.  The `left'
point is designated $x_1$. Let $N$ be the number of points of the Poisson point
process $P$ of intensity $\lam$ lying in $A^\prime$. Then the probability that $(m-l)$ sensors cover $x_1$ only
and $(n-l)$ cover $x_2$ only and $l$ sensors cover both is
$\prob{M_1 = l, M_2 = m-l,
  M_3 = n-l, M_4 = N-m-n+l}$ with
\[ (M_1,M_2,M_3,M_4) \sim \mbox{ Multinomial}(N , a_1 , a_2, a_3, a_4 = 1 - a_1 - a_2 - a_3),\]
where
\begin{eqnarray*}
a_i&  = & \frac{\EXP{\|B_i\|}}{\|A^\prime\| }, \qquad i=1,2,3.
\end{eqnarray*}
Since $N$ is Poisson with mean $\lam\|A^\prime\|$, the unconditional probability that $(m-l)$ sensors cover $x_1$ but not $x_2$, 
$(n-l)$ cover $x_2$ but not $x_1$, and $l$ sensors cover both $x_1$ and $x_2$ will be
\begin{eqnarray*}
  & &
  \EXP{\frac{N!}{l!(m-l)!(n-l)!(N-m-n+l)!}a_1^la_2^{(m-l)}a_3^{(n-l)}
    a_4^{(N-m-n+l)}} = \frac{ \left( \lambda \EXP{\|B_1\|} \right)^l}{l!}
  e^{-\lambda \EXP{\|B_1\|}} \\
  & & \times \frac{ \left(
      \lambda \EXP{\|B_2\|}\right)^{m-l}}{(m-l)!} e^{-\lambda\EXP{\|B_2\|}}
  \frac{\left( \lambda \EXP{\|B_3\|} \right)^{n-l}}{(n-l)!}e^{-\lambda
    \EXP{\|B_3\|}}.
\end{eqnarray*}
Hence from above calculations,
\begin{eqnarray}
  \EXP{\chi_m(x_1)\chi_n(x_2)}
  &=& \sum_{l=0}^{m\wedge n} \frac{ \left( \lambda \EXP{\|B_1\|}
    \right)^l}{l!} e^{-\lambda \EXP{\|B_1\|}} \times \frac{ \left(
      \lambda \EXP{\|B_2\|}\right)^{m-l}}{(m-l)!} e^{-\lambda
    \EXP{\|B_2\|}} \no  \\
  & & \hspace{0.3in} \times  \frac{\left( \lambda
      \EXP{\|B_3\|} \right)^{n-l}}{(n-l)!}e^{-\lambda
    \EXP{\|B_3\|}} \no \\
  & = & e^{-2\lambda \beta}e^{\lambda\EXP{\|B_1\|}}\sum_{l=0}^{m\wedge n}
  \frac{\left( \lambda \EXP{\|B_1\|} \right)^l}{l!} \times
  \frac{\left( \lambda \EXP{\|B_2\|}\right)^{m+n-2l}}{(m-l)!(n-l)!} .
  \label{eq:expexted-chi_m_n}
\end{eqnarray}
The last equality is obtained from the identities $\EXP{\|B_1\| +
  \|B_2\|} = \EXP{\|B_1\| + \|B_3\|} = \beta.$ We then have
\begin{eqnarray}
  \COV{ V_k(x_1) , V_k(x_2)} & = & \sum_{m,n = 0}^{k-1}
  \COV{\chi_m(x_1), \chi_n(x_2)} \nonumber \\
  & = & \sum_{m,n = 0}^{k-1} \left( \EXP{
      \chi_m(x_1 ) \chi_n(x_2)} - \EXP{\chi_m(x_1)} \EXP{\chi_n(x_2)}
  \right)  \no \\
  & = & e^{-2\lambda \beta} \sum_{m,n = 0}^{k-1} \left(
    e^{\lambda\EXP{\|B_1(x_1-x_2)\|}} \sum_{l=0}^{m \wedge n} \frac{ \left(
        \lambda  \EXP{\|B_1(x_1-x_2)\|} \right)^l}{l!} \right. \no \\
  & &  \hspace{0.5in} \left. \times \frac{ \left( \lambda E{\|B_2(x_1-x_2)\|}
      \right)^{m+n-2l}}{(m-l)!(n-l)!} -  \frac{ \left( \lambda
        \beta\right)^{m+n}}{m!n!} \right),
  \label{eq:COV-of-x1-x2}
  \end{eqnarray}
and
\begin{equation}
  \VAR{V_k} = \int_{R \times R} \COV{V_k(x_1), V_k(x_2)} \ dx_1 \ dx_2
  \label{eq:VAR-of-V_k}.
\end{equation}
\subsection{Limit Laws}
\label{sec:limit_laws_k_coverage}
We are now ready to obtain the
limit laws by letting $\lambda \rar \infty$ and scaling the shapes by $\delta.$ Let
$\mathcal{C}(\delta, \lambda)$ be the PB model $\mathcal{C}$ in which
the shapes are scaled by $\delta$, i.e., the shapes have the same
distribution as $\delta C.$ Let $V_k = V_k(\lambda,\delta)$ be the
resulting $k$-vacancy in $R.$
\begin{theorem}
  \label{thm:SL}
  If $\delta \to 0$ as $\lambda \to \infty $ such that $\delta^d
  \lambda \to \rho $ for $0 \leq \rho < \infty $ in the scaled coverage process $C(\del, \lam)$, then
\[
  V_k \rar \|R\|e^{-\rho \beta} \sum_{j=0}^{k-1} \frac{(\rho \beta
    )^j}{j!} \qquad \mbox{a.s.}
  \]
\end{theorem}
\begin{theorem}
  \label{thm:WL}
  Consider the scaled coverage process $\mathcal{C}(\delta,\lambda)$.
  If $\delta \to 0$ as $\lambda \to \infty $ such that $\delta^d
  \lambda \to \rho $ where $0 \leq \rho < \infty $, then
  \begin{eqnarray}
    \EXP{V_k}  & \to & \|R\| \sum_{m = 0}^{k-1}
    \frac{e^{-\rho \beta}(\rho \beta)^m}{m!},
    \label{eq:limit-expected-V_k}  \\
    \EXP{|V_k  - \EXP{V_k }|^p} & \to & 0  ~~~~~~~~~ \mbox{for
      $1 \leq p < \infty$}, \label{eq:limit-pth-moment} \\
    \lambda \VAR{V_k} & \to &  \sigma^2
    \label{eq:limit-V_k-variance},
  \end{eqnarray}
where
  \begin{eqnarray}
    \label{eq:sigma-square}
    \sigma^2 & = &\rho \|R\|e^{-2\rho \beta}\left(\sum_{m,n =
        0}^{k-1} \int_{B_{0,2\tau}}\left( e^{\rho \EXP{\|B_1(y)\|}}
        \sum_{l=0}^{m\wedge n} \frac{ \left( \rho \EXP{\|B_1(y)\|}
          \right)^l}{l!} \times \frac{\left( \rho
            \EXP{\|B_2(y)\|}\right)^{m+n-l}}{(m-l)!(n-l)!}
      \right. \right. \no \\
    & & \left. - \left. \frac{\left( \rho
            \beta\right)^{m+n}}{m!n!}\right) dy \right).
  \end{eqnarray}

\end{theorem}
%
%
\begin{theorem}
  \label{thm:CLT}
  If $\delta \to 0$ as $\lambda \to \infty $ such that $\delta^d
  \lambda \to \rho $ for $0 \leq \rho < \infty $ in the scaled coverage process $C(\del, \lam)$, then
  \[
  \sqrt{\lambda}\left( V_k-\EXP{V_k} \right) \rar N(0,\sg^2),
  \]
  in distribution where $\sg^2$ is as defined in
  (\ref{eq:sigma-square}).
\end{theorem}
It follows from Theorem~\ref{thm:SL} that if $\delta^d \lambda \rar
\infty$, then $V_k \rar 0$ almost surely. However, note that
Theorems~\ref{thm:SL} and ~\ref{thm:WL} do not guarantee complete
coverage with high probability for large enough $\lambda.$ We now
consider such a requirement, key to which is the inequality (\ref{coverage_bounds}) below.

For the following two theorems, we assume that $d = 2$, the operational area $R$ to
be the unit square $[0,1]^2,$ and the sensing areas to
be discs of radius $r>0$ satisfying $\pi r^2 \leq 1$. This last requirement is only to give a compact expression in the inequality below. Since our interest is in the asymptotic behavior of the
coverage process as $\lambda \to \infty$ and $r = r_{\lambda} \to 0,$ this is satisfied for all large enough $\lambda.$
For any $\lambda, r > 0$ define the event $Z_{\lambda}(r) := \{ V_k(\lambda, r) > 0 \} $.
For radii $r_{\lam} > 0$ which are decreasing in $\lambda$ to zero, we shall
abbreviate $Z_{\lambda}(r_{\lambda})$ by $Z_\lambda$. We have the following inequality.
\begin{equation}
  \frac{1}{1+\theta} \, \leq \, \prob{Z_\lambda(r)} \leq 2 e^{-\lambda
    \pi r^2} \left( 1 +  \lambda^2 \pi r^2(1 + \frac{2}{\lam \pi r})
\sum_{i=0}^{k-1} \frac{(\lambda \pi r^2)^i}{i!} \right)
  \label{coverage_bounds},
\end{equation}
where
\[
\theta = \frac{4 (k+1)!}{e^{-\lambda \pi r^2} \lambda (\lambda \pi
r^2)^k }.
\]
A formal proof of (\ref{coverage_bounds}) is given in the appendix.
The above inequality is an extension of Theorem~3.11 of \cite{Hall88} for
$k > 1$.  A similar inequality is derived in \cite{Zhang04} (proof of Theorem~1) under the condition that the operational area $R$ is a square of side length $\ell \to \infty$, with intensity $\lambda = \lambda(\ell) \to \infty$ (see (\ref{eq:k-coverage-zhang})) and the sensing area of the sensors is one,
i.e., $\pi r^2 = 1.$ The following result follows immediately from (\ref{coverage_bounds}).
\begin{theorem}
  \label{thm:full_cov} Suppose that $d=2,$ and let $r_\lambda^2 = \frac{\log \lambda + k \log
    \log \lambda + c(\lambda)}{\pi \lambda}.$ As $\lambda \rar
  \infty$, $\prob{Z_\lambda} \rar 0$ if $c(\lambda) \rar \infty.$
  Further, if $c(\lambda) \rar c \in \Re$, then $\prob{Z_\lambda} > c_1$
 for some constant $c_1 \in (0,1)$.
\end{theorem}
Theorem \ref{thm:full_cov} gives the critical radius required for
complete $k$-coverage with probability approaching $1$ in two
dimensions. We now show that by taking the radius to be a bit larger
than that obtained from Theorem \ref{thm:full_cov}, we can get a
stronger and more stable complete coverage regime. The following
discussion will become easier if we assume $\lambda=n$ where $n$ is an
integer.  To make the above notion of strong complete coverage more
precise, define the critical radius for complete coverage as
\begin{equation}
  r^\ast_n := \inf\{r_n >0: V_{k,s}( n, r_n)
  = 0 \}, \label{def_d_n}
\end{equation}
where $V_{k,s}$ is the $k$-vacancy in the unit square.
\begin{theorem}
  Let $d=2,$ and let $V_{k,s}$ be the vacancy in the unit square.  Let
  $r_n^\ast$ be as defined above. Then, almost surely,
  \begin{equation}
    \lim_{n \rar \infty} \frac{\pi n
      (r_n^\ast)^2}{\log n + k \log \log n} = 1.
    \label{strong_critical_radius}
  \end{equation}
  \label{thm:strong_critical_radius}
\end{theorem}
\textbf{Remark:} Let $0 < \epsilon < 1.$ The above result implies that
by taking the radius
\begin{equation}
  r_n^2 = (1+ \ep) \frac{\log n + k \log \log n}{\pi
    n},
\label{sup_critical_r_n}
\end{equation}
the unit square will be almost surely, completely $k$-covered for {\it
  all} $n$ large enough. Thus, if $n$ is large, by taking
the above $r_n$, which is eventually larger than the one given in
Theorem~\ref{thm:full_cov}, we can ensure a complete $k$-coverage
regime that will not see vacancies even if the number of sensors is
increased (with corresponding decrease in $r_n$). Further, the above
result gives a strong threshold in the sense that if
\begin{equation}
  r_n^2 = (1 - \ep) \frac{\log n + k \log \log n}{\pi
    n },
  \label{sub_critical_r_n}
\end{equation}
then the unit square will not be completely $k$-covered for {\it all}
large enough $n$, almost surely.
\subsection{Proofs}
\label{sec:proofs_k-coverage}
{\bf Proof of Theorem \ref{thm:SL}:} To obtain the required result,
first note that for the $k$-vacancy in a unit cube $D,$ $V_k(D),$ the
expectation is given by $\EXP{V_k(D)} = e^{-\rho \beta}
\sum_{j=0}^{k-1} \frac{(\rho \beta )^j}{j!}$. Now observe that the two
scaling regimes---(1) $R$ is fixed with $\lambda \delta^d \rar \rho$
and (2) $\delta = 1$, $R_l = lR$, and $l \rar \infty$ with $P$ being
a Poisson point process of intensity $\rho$ are
equivalent. See Section~3.4 of \cite{Hall88} for more discussion on
this. Theorem~\ref{thm:SL} follows using the same steps as in the proof of Theorem~3.6
in \cite{Hall88}. \qed

{\bf Proof of Theorem \ref{thm:WL}:} Since the sequence $\{V_k(\del,
\lambda)\}_{\del \geq 0, \lambda \geq 0}$ is bounded,
(\ref{eq:limit-expected-V_k}) and (\ref{eq:limit-pth-moment}) follow
from Theorem \ref{thm:SL}.

(\ref{eq:limit-V_k-variance}) follows from
(\ref{eq:COV-of-x1-x2}), (\ref{eq:VAR-of-V_k}) if $\lambda
\int_{R \times R} \COV{\chi_m(x_1), \chi_n(x_2)} dx_1 \, dx_2$
converges to the product of the integral on the r.h.s of
(\ref{eq:sigma-square}) and $ \left( \rho \| R \| e^{-2 \rho
    \beta}\right) .$ This is shown below.

Let $B_1^{\delta}(y) = (y+\delta C) \cap \delta C$ and $B_2^{\delta}(y) = (y+\delta C) \cap \delta C^c$.
Making first the change of variable $x=x_1$ and $y=x_1-x_2$ in (\ref{eq:COV-of-x1-x2}) and then $y$ to $\del y$, we get
{\small
  \begin{eqnarray}
\lefteqn{
\lambda \int_R \int_{x-R} \COV{\chi_m(x), \chi_n(x_2)} dx_2 \,
      dx} \no \\
 & = &
\lambda \int_R \int_{x-R}e^{-2\lambda \del^d \beta}\sum_{m,n =
      0}^{k-1} \left( e^{\lambda\EXP{\|B^{\delta}_1(y)\|}} \sum_{l=0}^{m \wedge n}
      \frac{ \left( \lambda  \EXP{\|B^{\delta}_1(y)\|} \right)^l}{l!}
      \times \right. \no \\
&&  \hspace{0.5in} \left. \frac{ \left( \lambda
      \EXP{\|B^{\delta}_2(y)\|} \right)^{m+n-2l}}{(m-l)!(n-l)!} - \frac{
      \left( \lambda \beta\right)^{m+n}}{m!n!} \right) dy \ dx \no \\
&=&\lambda \int_R \int_{x-R}  e^{-2\lambda \del^d \beta} \sum_{m,n =
      0}^{k-1} \left( e^{\lambda \del^d \EXP{\|B_1(\del^{-1}y)\|}} \sum_{l=0}^{m \wedge n}
      \frac{ \left( \lambda \del^d \EXP{\|B_1(\del^{-1}y)\|} \right)^l}{l!}
      \times \right. \no \\
&&  \hspace{0.5in} \left. \frac{ \left( \lambda \del^d
\EXP{\|B_2(\del^{-1}y)\|}
        \right)^{m+n-2l}}{(m-l)!(n-l)!} -  \frac{ \left( \lambda
         \del^d \beta\right)^{m+n}}{m!n!} \right) dy \ dx \no 
\end{eqnarray}
\begin{eqnarray}
    & = & \lambda \del^d e^{-2 \delta^d \lambda \beta}  \int_{R} \
    \int_{\del^{-1}(x-R)}\left( e^{\del^d  \lambda \EXP{\|B_1(y)\|}}
      \sum_{l=0}^{m\wedge n}\frac{\left(  \del^d\lambda
          \EXP{\|B_1(y)\|}\right)^l}{l!} \times\frac{\left( \del^d
          \lambda \EXP{\|B_2(y)\|}\right)^{m+n-2l}}{(m-l)!(n-l)!}
    \right.\no\\
    &&\hspace{0.5in}
    \left. - \frac{\left( \del^d \lambda \beta\right)^{m+n}}{m!n!}
    \right) dy  \ \ dx \no \\
    & = & \lambda \del^d e^{-2 \delta^d \lambda \beta} \int_{R} \
    f_\delta(x) \ dx, \label{eqn:cov_int}
  \end{eqnarray}
}
\normalsize
where
\begin{eqnarray}
  f_\delta(x) &:=&  \int_{\del^{-1}(x-R) \cap B_{0,2 \tau}}\left( e^{\del^d
      \lambda \EXP{\|B_1(y)\|}} \sum_{l=0}^{m\wedge n}\frac{\left(
        \del^d\lambda \EXP{\|B_1(y)\|}\right)^l}{l!}
    \times\frac{\left( \del^d \lambda
        \EXP{\|B_2(y)\|}\right)^{m+n-2l}}{(m-l)!(n-l)!} \right.\no\\
  &&\hspace{0.5in}
  \left. - \frac{\left( \del^d \lambda \beta\right)^{m+n}}{m!n!}
  \right) dy.
  \label{eqn:f-delta-of-x}
\end{eqnarray}
$f_{\del}(x)$ as well as the integrand in (\ref{eqn:f-delta-of-x}) are uniformly bounded, since for any $\ep > 0$ and all $\lam > 0$ sufficiently large, we have $\del^{-1}(x-R) \cap B_{0,2 \tau} \subset B_{0,2 \tau}$, $|\del^d \lam - \rho| \leq  \epsilon$
and $B_i \subset B_{0,\tau}.$  Therefore, from (\ref{eqn:cov_int}), (\ref{eqn:f-delta-of-x}) and the dominated convergence theorem,
we have
\small{
  \begin{eqnarray}
    f_\delta(x) & \to & c_0 :=  \int_{B_{0,2\tau}}\left( e^{\rho
        \EXP{\|B_1(y)\|}} \sum_{l=0}^{m\wedge n} \frac{\left(
          \rho \EXP{\|B_1(y)\|}\right)^l}{l!} \times \frac{\left( \rho
          \EXP{\|B_2(y)\|}\right)^{m+n-2l}}{(m-l)!(n-l)!}
      -\frac{\left( \rho \beta \right)^{m+n}}{m!n!}\right) dy,
    \nonumber \\
    \label{conv_f_delta}
  \end{eqnarray}
}
\normalsize 
and
\[
 \lambda  \int_{R \times
  R}\COV{\chi_m(x_1)\chi_n(x_2)} \ dx_1 \ dx_2 \rar \rho e^{-2 \rho \beta} \int_R c_0 \
dx = \rho e^{-2 \rho \beta} c_0 \|R\|,
\]
as $\lambda \rar \infty.$ This proves
(\ref{eq:limit-V_k-variance}). \qed \\
\begin{figure}[tb]
  \centerline{\includegraphics[scale=0.7]{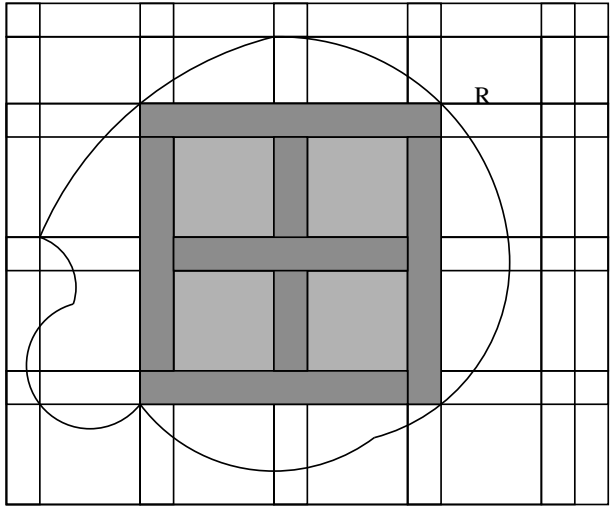}}
  \caption{Figure shows the sets $A_1$, $A_2$ and $A_3.$ The squares
    are shaded light gray. The squares are of length $\tau r \delta$ and the
    `spacings' between the sqaures are of width $2 \tau \delta$.  The
    light gray shaded region within $R$ is $A_1$, the darker region within
    $R$ is $A_2$ and the unshaded region within $R$ is $A_3$
    (Adapted from \cite{Hall88}). }
  \label{fig:CLT-illustrations}
\end{figure}

{\bf Proof of Theorem \ref{thm:CLT}:} Let $r$ be a large positive
constant. Divide all of $\Re^d$ into a regular lattice of
$d$-dimensional cubes of side length $(\tau r \delta),$ with each cube
separated from its adjacent cubes by a $d$-dimensional `spacing' of
width $2 \tau \delta.$ Let $A_1$ denote the union of those cubes which
are wholly within $R,$ $A_2$ the union of the spacings that are wholly
contained in $R$ and $A_3$ the intersection of $R$ with all those
cubes and spacings that are only partially within $R$.
Fig.~\ref{fig:CLT-illustrations} illustrates these sets for $d=2.$
Since $A_1,$ $A_2$ and $A_3$ form a partition of $R$, the vacancy
$V_k$ within $R$ may be written as,
\[
V_k = V^{(1)}_k + V^{(2)}_k + V^{(3)}_k,
\]
where $V^{(i)}_k$ is the $k$-vacancy within the region $A_i.$ Under
the assumptions of the theorem, as $\delta \to 0,$ the cubes get
finer. Further the number of spacings is less than $\| R \|/(r \tau
\delta)^d.$ Since the volume of the `spacings' is $(\tau r
\delta)^{d-1} \times (2 \tau \delta),$ we have
\begin{equation}
  \label{12}
  \|A_3 \| \rightarrow 0, \qquad \mbox{ and } \qquad \|A_2\| \leq
  2 \| R \| r^{-1}.
\end{equation}
From (\ref{eq:VAR-of-V_k}), we get
\begin{eqnarray}
  \VAR{V^{(i)}_k} & = & \sum_{m,n =0}^{k-1}
  \int_{A_i}dx_1
  \int_{A_i} \COV{\chi_m(x_1)\chi_n(x_2)}dx_2 \no \\
  & = & \delta^d \sum_{m,n =0}^{k-1} \int_{A_i} e^{-2 \delta^d
    \lambda \beta} f_{\delta}(x_1)dx_1
  \leq  M \delta^d k^2 \| A_i \|,
  \label{eq:V_(i)-variance-bounds}
\end{eqnarray}
where the last inequality follows from (\ref{conv_f_delta}).  From the
inequality above and (\ref{12}), we get
\begin{equation}
  \label{13}
  \lim_{\lambda \rar \infty} \lambda \VAR{V^{(3)}_k} = 0.
\end{equation}
From (\ref{12}) and (\ref{eq:V_(i)-variance-bounds}), we also get
\begin{equation}
  \label{14}
  \lim_{r \rar \infty} \limsup_{\lambda \rar \infty} \lambda
  \VAR{V^{(2)}_k} = 0.
\end{equation}
Thus, for any $\epsilon >0,$ we can choose $r$ large enough so
that $\lambda \VAR{V^{(2)}_k} < \epsilon,$ for all $\lambda$
sufficiently large. Hence in order to obtain the central limit theorem
that we are seeking, we need to concentrate on $V^{(1)}_k$ and obtain a
central limit theorem for it. Since $\lambda \VAR{V_k} \rar \sigma^2$
from (\ref{eq:limit-V_k-variance}), we need to show the following.
\begin{equation}
  \label{15}
  \frac{  V^{(1)}_k - \EXP{V^{(1)}_k} }{ \left( \VAR{V^{(1)}_k}
    \right)^{1/2}} \stackrel{d}{\rar} N(0,1),
  \qquad  \lim_{r \rar \infty}\limsup_{\lambda \rar
    \infty}|\lambda (\VAR{V^{(1)}_k} - \VAR{V_k})| = 0.
\end{equation}
Let $n = n(\lambda)$ denote the number of cubes of side-length $\tau r
\delta$ in $A_1$, and let $D_i$ denote the $i$-th such cube, $1 \leq i
\leq n.$ Denoting the $k$-vacancy in $D_i$ by $U_i$, we have $V^{(1)}_k
= \sum_i U_i$. Under the scaling regime, each shape is contained in
$B_{0, \tau \delta}$ and the spacing between the cubes is $2 \tau
\delta$, and so no shape can intersect more than one cube.  Hence,
given $\lambda,$ the $U_i'$s are independently distributed and we have
\begin{displaymath}
  \VAR{V^{(1)}_k} = \sum_i \VAR{U_i} = n\VAR{U_1}. \no
\end{displaymath}
Let $D$ be a $d$-dimensional cube of side $\tau r$ with the same
orientation as $D_1$. For any two real sequences $a_n,b_n$ $a_n \sim
b_n$ implies that $a_n/b_n \to 1$ as $ n \to \infty.$ From
(\ref{eq:COV-of-x1-x2}), (\ref{eq:VAR-of-V_k}) we can write
\begin{eqnarray}
  \VAR{V^{(1)}_k} & = & n e^{-2\delta^d \lambda \beta} \sum_{m,n =
    0}^{k-1}\int \int_{D_1 \times D_1} \left( e^{ \lambda
      \EXP{\|(x_1-x_2+\delta S) \cap \delta S\|}}
    \sum_{l=0}^{m \wedge n}\frac{ \left( \lambda
        \EXP{\|(x_1-x_2+\delta S) \cap \delta S\|}\right)^l}{l!}
  \right. \no \\
  & & \left. \times \frac{\left( \lambda
        \EXP{\|(x_1-x_2+\delta S) \cap (\delta
          S)^c\|}\right)^{m+n-2l}}{(m-l)!(n-l)!}
    - \frac{\left( \delta^d \lambda \beta\right)^{m+n}}{m!n!}
  \right) dx_1 dx_2  \no \\
  & \sim & n\delta^{2d} e^{-2\rho \beta} \sum_{m,n =
    0}^{k-1}\int \int_{D \times D} \left( e^{\rho \EXP{\|B_1\|}}
    \sum_{l=0}^{m \wedge n}\frac{ \left( \rho
        \EXP{\|B_1\|}\right)^l}{l!} \times \frac{\left( \rho
        \EXP{\|B_2\|}\right)^{m+n-2l}}{(m-l)!(n-l)!}  \right. \no \\
  & & \left. - \frac{\left( \rho \beta\right)^{m+n}}{m!n!} \right)
  dx_1 dx_2.  \label{order_of_var_V_1}
\end{eqnarray}
The rest of the proof follows as in (\cite{Hall88} pp.  157-158). \qed

{\bf Proof of Theorem \ref{thm:strong_critical_radius}:} Recall that
$Z_n(r_n)= \{V_k(n,r_n) > 0 \}.$ Suppose we show that for the choice
of $r_n$ as in (\ref{sup_critical_r_n}), we have $\prob{Z_n(r_n)
  \mbox{ infinitely often }} = 0,$ then it follows that
\begin{equation}
  \limsup_{n \rar \infty} \frac{\pi n (r_n^\ast)^2}{\log n + k \log \log n}
  \leq 1. \label{lim_sup}
\end{equation}
On the other hand, for the choice of $r_n$ as in
(\ref{sub_critical_r_n}), if we show that $\prob{Z_n^c(r_n) \mbox{
    infinitely often }} = 0,$ then we can conclude that
\begin{equation}
  \liminf_{n \rar \infty} \frac{\pi n (r_n^\ast)^2}{\log n + k \log \log n}
  \geq 1. \label{lim_inf}
\end{equation}
First we show (\ref{lim_sup}). Take subsequence $n_j = j^a,$ where
$a > 0$ will be chosen appropriately later. Define the events $E_j =
\cup_{n = n_j}^{n_{j+1}} Z_n(r_n) \subset Z_{n_{j}}(r_{n_{j+1}}).$ Then,
$\prob{E_j}  \leq  \prob{Z_{n_{j}}(r_{n_{j+1}})}$.

For $r_n$ as in (\ref{sup_critical_r_n}), observe that $n_j
r_{n_{j+1}}^2 \rar \infty.$ Hence, from (\ref{coverage_bounds}), we
get, for large enough $j$,
\[
\prob{Z_{n_{j}}(r_{n_{j+1}})} \leq \zeta_1 k e^{-n_j \pi
  r_{n_{j+1}}^2} n_j (n_j \pi r_{n_{j+1}}^2)^k,
\]
where $\zeta_1$ is some constant. Choose $\eta > (1+\ep)^{-1}.$ Since
$n_j/n_{j+1} \uparrow 1,$ as $j \rar \infty,$ we have for all large
enough $j$, $n_j/n_{j+1} \geq \eta.$ Hence, $\eta (1+ \ep) \log
n_{j+1} \leq n_j \pi r_{n_{j+1}}^2 \leq 4 \log n_{j+1},$ so that
\[
\prob{Z_{n_{j}}(r_{n_{j+1}})} \leq \zeta_2 e^{-\eta(1+\ep) \log
  n_{j+1}} (\log n_{j+1})^k n_j \leq \zeta_3 (\log j)^k (j+1)^{-a(\eta
  (1+\ep) - 1)},
\]
where $\zeta_2, \zeta_3$ are constants that depend on $a,k,\ep.$
Thus for $a$ large enough we get $\sum_{j=1}^{\infty} \prob{E_j} <
\infty.$ By the Borel-Cantelli Lemma, $E_j$ happens only finitely
often almost surely. By the definition of $E_j$, this implies that
$\prob{Z_n(r_n) \mbox{ infinitely often }} = 0.$

Similarly using the lower bound for $\prob{Z_n}$, we can prove
(\ref{lim_inf}). As above, let $n_j = j^a$, $j \geq 1$ and $a >0$ to be
chosen later. Let $r_n$ be as in (\ref{sub_critical_r_n}).  Using the
lower bound in (\ref{coverage_bounds}), we get
\[
\prob{\cup_{n = n_j}^{n_{j+1}} Z_n^c(r_n)} \leq
\prob{Z_{n_{j+1}}^c(r_{n_j})} \leq \frac{\theta_j}{1+\theta_j} \leq
\theta_j,
\]
where $$ \theta_j = \frac{4(k+1)!}{e^{-n_{j+1}\pi
    r_{n_j}^2}n_{j+1}(n_{j+1}\pi r_{n_j}^2)^k}.$$ Proceeding as in
the proof of the upper bound, we can show that the above probability
is summable. This completes the proof of the theorem.
\hfill \qed
\section{Path Coverage in the Markov-Poisson-Boolean Model}
\label{sec:markov-boolean}
\subsection{Preliminaries}
\label{sec:prelims_2}
Let $\mathcal{C}_t= \{Y_i^t(X_i + C_i), i \geq 1 \}$ be the MPB model
as defined in Section \ref{sec:Markov-Boolean}. Recall that $\{Y^t_i, t \geq 0 \}$
are i.i.d. copies of
$\{Y^t, t \geq 0 \}$ which is a $\{0,1\}$-valued continuous time stationary
Markov process. The sets $C_i$ are independent and distributed as $C$ satisfying
$C \subset B_{0,\tau}$ for some fixed $\tau >0$.
Let $\mu_0$ be the transition rate from the $0$-state
to the $1$-state and $\mu_1$ the transition rate from $1$-state to the
$0$-state.  This of course means that in each visit, $Y^t$ is in the
$0$ and $1$ states for exponentially distributed times with parameters
$\mu_0$ and $\mu_1$ respectively. Then, (see \cite{Ross05}, Chapter
6), the stationary probability of the sensor being in state $j$ is
given by $p_j := \prob{Y^t = j} = \frac{\mu_{1-j}}{\mu_1 + \mu_0}$,
$j \in \{0,1\}.$ The transition probabilities between the states are
defined by $p_t(j,k) := \prob{Y^{s+t} = k|Y^s = j}$ for $j,k \in
\{0,1\}.$ It can be shown that
\begin{displaymath}
  p_t(j,j) = (1-p_j) e^{-\gamma t} + p_j,
\end{displaymath}
where $ \gamma = \mu_0 + \mu_1$ and $p_t(j,k) = 1-p_t(j,j)$ for $j,k
\in \{0,1\}$ and $j \neq k.$

Now, consider a target moving on a straight~line path in $\Re^d$ of
length $L$ units with a velocity of $c$ units per second.  As in the
previous section, we will let $X_i=(X^1_i,X^2_i, \ldots ,X^d_i)$ be a
point in $\Re^d.$ Without loss of generality, we can consider the target to be moving
along a coordinate axis, say $X^1,$ for $T:=L/c$ units of time.
Define $\us := (cs,0,\ldots,0) \in \Re^2,$ where $s \in \Re^+.$ $\us$
is the position of the target on the $X^1$-axis at time $s.$ Let
$\zeta(s)$ be the indicator variable for the target being covered by
one or more of the sets $C_i$ at time $s,$ and $V_T= V_{1,T}$ be the random
variable denoting the duration for which the target is not covered in
$[0,T]$, i.e.,
\begin{equation}
  V_T := \int_0^T (1-\zeta(s)) \ ds.
  \label{eqn:V_T}
\end{equation}
The effect of the transitions of $Y^t_i$ makes the study of $V_T$
interesting. We first calculate the expectation and variance of $V_T$.
Since the $Y^t_i$ are stationary, the expectation of $V_T$ is
straightforward and is given by
\begin{eqnarray}
  \EXP{V_T} = \int_0^T\prob{\us \notin \cC_s}ds = Te^{-\lambda p_1
  \beta},
  \label{eqn:expected-V_T}
\end{eqnarray}
where $\beta = \EXP{\|C\|}.$ The second moment of $V_T$ can be
obtained from
\begin{eqnarray*}
  \EXP{V_T^2} &=&  \int_0^T \int_{0}^T \prob{\us_1 \notin
    \cC_{s_1}, \us_2 \notin \cC_{s_2}}\ ds_1 \ ds_2 \\
  & = & 2 \int_0^T \int_{s_1}^T \prob{\us_1 \notin
    \cC_{s_1}, \us_2 \notin \cC_{s_2}}\ ds_1 \ ds_2.\\
\end{eqnarray*}
We now evaluate this integral. Define $B_1(\us_1,\us_2) := (\us_1 + C) \cap (\us_2
+ C),$ $B_2(\us_1,\us_2) := (\us_1 + C) \cap (\us_2 + C)^c$ and $B_3(\us_1,\us_2) := (\us_2 + C)
\cap (\us_1 + C)^c.$ From this definition, $B_1,$ $B_2$ and $B_3$ are
disjoint and we can see that $\|B_2\|$ and $\|B_3\|$ have the same
distribution and $\|B_1\| = \|(\us_1 - \us_2+ C)\cap C\|$. We will suppress the
obvious arguments of the sets $B_i$. As in the
arguments leading to (\ref{eq:expexted-chi_m_n}), we consider a set $A^\prime =
B_{\us_1,2\tau} \cup B_{\us_2,2\tau}$ and assume there are $N$ sensors
within $A^\prime$.  Then the integrand in the above equation is
\[ \mathsf{Pr} \left( \mbox{0 on sensors covering $\us_1$ at
 $s_1$; 0 on sensors covering $\us_2$ at $s_2$} \right), \]
which equals $\prob{M_1 = 0, M_2 = 0, M_3 = 0,M_4 = N}$, where $M_1 =$
no: of sensors active at $s_1$ or $s_2$ among those in $A^\prime$
which can cover $\us_1$ and $\us_2$, $M_2 =$ no: of sensors active at
$s_1$ among those in $A^\prime$ which can cover $\us_1$ and not
$\us_2$, $M_3 =$ no: of sensors active at $s_2$ among those in
$A^\prime$ which can cover $\us_2$ and not $\us_1$, $M_4 =$ remaining
sensors in $A^\prime$. Now note that
\[ (M_1,M_2,M_3,M_4) \sim \mbox{ Multinomial}(N , a_1 , a_2, a_3, a_4 = 1 - a_1 - a_2 - a_3),\]
where
\begin{eqnarray*}
a_1&  = & (1 - p_0p_{s_2-s_1}(0,0))\EXP{\|B_1\|}/\|A^\prime\|, \\
a_2 & = & p_1\EXP{\|B_2\|}/\|A^\prime\|, \\
a_3 & = & p_1\EXP{\|B_3\|}/\|A^\prime\|.
\end{eqnarray*}
 $a_1$ is the probability that a sensor
can cover both of $\us_1$ and $\us_2$ and is active at either $s_1$ or
$s_2$. It is the product of probabilities of a sensor lying in $B_1$
and making the desired transition of states. The first probability is
$\EXP{\|B_1\|}/ \|A^\prime\|$. The second probability is obtained as
follows:
\begin{eqnarray*}
\prob{\mbox{a given sensor is on at $s_1$ or $s_2$}} & = &
\prob{\mbox{the sensor is on at $s_1$}} \\
& & + \prob{\mbox{the sensor is off at $s_1$ and on at $s_2$}} \\
& = & p_1 + p_0p_{s_2-s_1}(0,1) = 1 - p_0p_{s_2-s_1}(0,0).
\end{eqnarray*}
The remaining $a_i$'s are even simpler to obtain by similar
calculations. Since $N$ follows a Poisson distribution with mean $\lam \|A^\prime \|$,
the integrand in the above expression for $ \EXP{V_T^2}$ is given by
\[  \sum_{\ell = 0}^{\infty} a_4^{\ell} \frac{(\lam \| A^\prime \|)^{\ell}}{\ell !} e^{-\lam \| A^\prime \|}
= e^{-\lam(a_1 + a_2 + a_3) \| A^\prime \| } .\]
\begin{eqnarray*}
  \EXP{V_T^2} &=& 2\int_0^T \int_{s_1}^T \left(e^{-\lambda
      p_1\EXP{\|B_2\|}} \right) \left(e^{-\lambda
      p_1 \EXP{\|B_3\|}} \right) \left( e^{-\lambda
      \EXP{\|B_1\|}(1-p_0p_{s_2-s_1}(0,0))} \right)  ds_2 \; ds_1 \\
  &=& 2\int_0^T \int_{s_1}^T e^{-2p_1 \lambda \beta }  e^{-\lambda
    \EXP{\|B_1\|}(1-2p_1-p_0p_{s_2-s_1}(0,0))} ds_2 \; ds_1.
\end{eqnarray*}
Hence,
\begin{eqnarray}
  \VAR{V_T} & = &  \EXP{V_T^2} - (\EXP{V_T})^2 \no \\
  & =  & 2e^{-2 p_1 \lambda \beta} \int_0^T \int_{s_1}^T \left(
    e^{-\lambda \EXP{\|B_1\|}(1-p_0p_{s_2-s_1}(0,0)-2p_1)} - 1
  \right)ds_2 \; ds_1
  \label{eqn:var(V_T)}.
\end{eqnarray}
Observe from (\ref{eqn:expected-V_T}) that the expected coverage
depends neither on the shapes of the sensing regions nor on the
transition probabilities of the on-off process. Their effects
are picked up only in the variance of the coverage as can be seen in
(\ref{eqn:var(V_T)}).
\subsection{Limit Laws}
\label{sec:limit_laws_2}
We can now delineate the scaling laws. Let $\cC_t(\delta,\lambda)$
be the process $\cC_t$ in which the shapes $C$ are scaled by
$\delta.$ In addition, the parameters of the Markov chain governing
the on-off process will also be made to depend on $\delta$. Without
loss of generality, let $ \mu_0(\delta), \mu_1(\delta)$ be such that
$\mu_0(\delta) / \mu_1(\delta) \equiv \mu_0/\mu_1$ for all $\delta.$
It can be easily seen that the proofs remain valid when
$\mu_0(\delta) / \mu_1(\delta) \rar \mu_0/\mu_1$ as $\del \rar 0.$
Hence the scaled stationary probability is $p_i(\delta) \equiv p_i$
and the scaled transition probabilities are $p_s^{\del}(0,0) = p_1
e^{-\gamma(\delta)s} + p_0,$ where $\gamma(\delta) = \mu_0(\delta) +
\mu_1(\delta).$
Observe that the case
$\mu_0(\delta) / \mu_1(\delta) \rar 0 (\infty)$ imply that $p_0 = 1
(0)$ and they correspond to sensors being always off(on).

We will first state the four main theorems for the MPB model and
provide their proofs after a brief discussion.

\begin{theorem}
  \label{thm:scal_var}
  Consider the scaled coverage process
  $\mathcal{C}_t(\delta,\lambda).$ Let $\delta \to 0$ as $\lambda \to
  \infty $ such that $\delta^d \lambda \to \rho $ where $0 \leq \rho <
  \infty $ and $\del \gamma(\delta) \to a_0$, where $0 \leq a_0 \leq
  \infty.$ Then,
  \begin{eqnarray}
    \label{eqn:lt_exp}
    \EXP{V_T}& \to & Te^{-p_1 \rho \beta}, \\
    \label{eqn:lt_L_p}
    \EXP{|V_T - \EXP{V_T}|^\mathrm{p}} & \to & 0 \qquad \mbox{for $\;\; 1
      \leq \mathrm{p} < \infty.$}
  \end{eqnarray}
  \begin{equation}
    \label{eqn:lt1_var}
   \delta^{-1} \VAR{V_T(\delta)}   \to  \sigma_1^2(a_0),
  \end{equation}
where
  \begin{equation}
   \label{eqn:sg1_a0}
    \sigma_1^2(a_0) = 2Te^{-2 p_1 \rho \beta} \int_0^{\frac{2\tau}{c}} \left(
      e^{p_1 \rho \EXP{\|(\us+C)\cap C\|}(p_1 + p_0e^{-a_0s})} - 1 \right) \ ds.
  \end{equation}
\end{theorem}
\begin{theorem}
  \label{thm:unrel_clt} Let the assumptions be as in Theorem
  \ref{thm:scal_var}.
Let $0 \leq a_0 \leq \infty$. Then,
    \begin{equation}
      \label{eqn:clt1}
      \delta^{- \frac{1}{2}}(V_T(\delta) - \EXP{V_T(\delta)}) \rar  N(0,\sigma_1^2),
    \end{equation}
in distribution, as $\delta \to 0.$
%
    %
    %
\end{theorem}
The following two theorems describe the asymptotic behavior of the coverage process
$V_T$ as the size of the operational area becomes very large.
\begin{theorem}
  \label{thm:unrel_sl}
 Consider the coverage process $\mathcal{C}_t = \mathcal{C}_t(1,\lambda)$. As $T \to \infty,$
  \begin{displaymath}
    \frac{V_T}{T} \to \frac{e^{-p_1 \lambda \beta}}{c} \qquad \mbox{a.s.}
  \end{displaymath}
\end{theorem}
\begin{theorem}
\label{thm:unrel_clt2} For the coverage $\mathcal{C}_t = \mathcal{C}_t(1,\lambda)$, as $T \rar \infty$,
\begin{equation}
\label{eqn:unrel_clt2} T^{-1/2}(V_T - \EXP{V_T}) \rar
N(0,\sigma_2^2),
\end{equation}
in distribution, where
\begin{equation*}
    \sigma_2^2 = 2e^{-2 p_1 \lam \beta} \int_0^{\frac{2\tau}{c}} \left(
      e^{p_1 \lam \EXP{\|(\us+C)\cap C\|}(p_1 + p_0e^{-(\mu_0 +
          \mu_1)s})} - 1 \right) \ ds.
\end{equation*}
\end{theorem}
%
%
\textbf{Remark:} The scaling of the shapes by $\del$ can also be
viewed as follows: If the velocity of the target increases by
$1/\del$, then the time it shall spend in the region of a sensor shrinks
by a factor $\del$ but it shall travel a larger distance and hence be
seen by more sensors. In light of this remark, one can view the
scaling results in the unreliable sensor networks case as tracking of
a high-velocity target in a highly fluctuating sensor network. Observe that
\begin{displaymath}
  \sigma_1^2(0) = 2Te^{-2p_1 \rho
    \beta}\int_0^{2\tau/c}\left(e^{p_1 \rho \EXP{\|(\us+C)\cap C\|}} -
    1 \right) \ ds.
\end{displaymath}
Since the target is moving fast, the contribution to the variance is coming from locations that
are close to each other. When $a_0$ is small, the target moves much faster as compared to
the rates of transition of the sensors between the on and the off
states. Thus given that there was a sensor on at a given location, the chance of it being on in the near future
will be close to 1. Similarly if $a_0 = \infty$, then
\begin{displaymath}
  \sigma_1^2(\infty) = 2Te^{-2p_1 \rho
    \beta}\int_0^{2\tau/c}\left(e^{p_1^2 \rho \EXP{\|(\us+C)\cap C\|}} -
    1 \right) \ ds.
\end{displaymath}
Hence when $a_0$ is large, the transitions happen at a much faster rate compared to the speed of the target. So
given that there was a sensor on at a location, the chance of the target seeing a sensor on in the near future will be roughly
the stationary probability of a sensor being on, which is $p_1.$
\subsection{Proofs}
\label{sec:MPB_proofs}
{\bf Proof of Theorem \ref{thm:scal_var}:} (\ref{eqn:lt_exp}) is
straightforward. For (\ref{eqn:lt_L_p}), as in the proof of
(\ref{eq:limit-pth-moment}), it suffices to show that the variance
converges to $0$. Putting $s = s_2 - s_1$ in (\ref{eqn:var(V_T)}) and
the expression for $p_{s_2-s_1}(0,0)$ and simplifying, we get
\begin{eqnarray*}
  \VAR{V_T(\delta)} = 2 e^{-2\delta^d\lambda p_1 \beta}\int_0^T
  \int_0^{T-s_1} \left(e^{\delta^d\lambda p_1\EXP{\|(\delta^{-1}\us +
        C) \cap C \|}(p_1 + p_0e^{-\gamma(\delta)s})} - 1 \right) \ ds \
        ds_1.
\end{eqnarray*}
%
Replacing $s$ by $s\del$ we get
%
\begin{equation}
  \label{eqn_scal_var1}
  \VAR{V_T(\delta)} =  2 \delta
  e^{-2 \delta^d \lambda p_1 \beta}\int_0^T
  \int_0^{\delta^{-1}(T-s_1)} \left( e^{\delta^d\lambda p_1
      \EXP{\|(\us+C)\cap C \|}(p_1 + p_0e^{-\delta \gamma(\delta)s})}
    - 1 \right) \ ds \  ds_1.
\end{equation}
\remove{
\begin{equation}
  \label{eqn_scal_var2}
   =   \frac{2}{\gamma(\delta)} e^{-2 \delta^d \lambda p_1 \beta}
  \int_0^T \int_0^{\gamma(\delta)(T-s_1)} \left( e^{\delta^d\lambda
      p_1 \EXP{\|(\frac{\us}{\del \gamma(\delta)}+C)\cap C \|}(p_1 +
      p_0e^{-s})} - 1 \right)\ ds \ ds_1.
\end{equation}}
The range of integration in the inner integral diverges
to $\infty$ as $\delta \to 0,$ and the limit of the integrand is as described in
(\ref{eqn:sg1_a0}) (using the fact that $C \subset B_{0,\tau}$).
We can justify the convergence to limits since the inner integral and the integrand are
uniformly bounded and hence the dominated convergence theorem applies. Indeed, for any $\epsilon >0$ fixed and
sufficiently large $\lam$ we have $\lam \del^d \leq \rho + \ep$ and
\begin{eqnarray*}
\int_0^{\delta^{-1}(T-s_1)} \left( e^{\delta^d\lambda p_1
      \EXP{\|(\us+C)\cap C \|}(p_1 + p_0e^{-\delta \gamma(\delta)s})}
    - 1 \right) \ ds & \leq & \int_0^{\delta^{-1}T} \left( e^{(\rho + \ep) p_1
      \EXP{\|(\us+C)\cap C \|}(p_1 + p_0)}
    - 1 \right) \ ds \\
 & \leq & \frac{2 \tau}{c} \left( e^{(\rho + \ep)p_1
      \EXP{\|C \|} }
    + 1 \right)  < \infty.
\end{eqnarray*}
\qed

{\bf Proof of Theorem \ref{thm:unrel_clt}:} We shall first prove
(\ref{eqn:clt1}).
Proof follows the same idea used for proving the central limit theorem
in case of $k$-coverage.

Choose a large constant $r$. Divide $[0,T]$ into alternating
intervals of lengths $r\del \tau$ (type~1) and $2\del \tau/c$
(type~2), where $c$ is the velocity of the target. Truncate the
interval containing $T$ at $T$. Denote the union of type~1 intervals
by $A_1$ and the union of type~2 by $A_2$. Let the vacancies arising
in $A_1$ be denoted by $V^{(1)}$ and that in $A_2$ by $V^{(2)}.$ Now, $V_T =
V^{(1)} + V^{(2)}.$ Note that $\|A_2\| \leq \left( 2 \delta \frac{\tau}{c}
\times
  \frac{T}{r \delta \tau} \right) = \frac{2T}{cr} \to 0$ as $r \to
\infty.$ In (\ref{eqn_scal_var1}), the inner integral converges and
hence it is bounded. Using this fact when the range of integration is
$A_i$, we get
\begin{displaymath}
  \VAR{V^{(i)}} \leq C\del\|A_i\|.
\end{displaymath}
It follows that
\begin{equation}
  \label{eqn:V_2_to_0}
  \lim_{r \to \infty}\limsup_{\delta \to 0} \ \frac{\VAR{V^{(2)}}}{\del}
  = 0.
\end{equation}
Thus we need to show the following to prove the central limit
theorem.
\begin{equation}
  \label{eqn:for_unrel_clt}
  \frac{  V^{(1)} - \EXP{V^{(1)}} }{ \left( \VAR{V^{(1)}}
    \right)^{1/2}} \stackrel{d}{\rar} N(0,1),
  \qquad  \lim_{r \rar \infty}\limsup_{\lambda \rar
    \infty}|\lambda (\VAR{V^{(1)}} - \VAR{V_T})| \to 0.
\end{equation}
Let $n = n(\delta)$ denote the number of intervals in $A_1$.  Let
$D_i$ denote the $i$-th interval of type~1.  Let $U_i$ be the vacancy
in $D_i.$ We can see that the $U_i$ are i.i.d. random variables except
for last one corresponding to the truncated interval and $V^{(1)} =
\sum_{i=1}^n U_i.$ As in (\ref{eqn_scal_var1}), we obtain,
\begin{eqnarray*}
  \VAR{V^{(1)}} & \sim & n\VAR{U_1} \no \\
  & = & n 2 \delta e^{-2\delta^d\lambda p_1\beta}\int_0^{r \delta \tau}
  \int_0^{\delta^{-1}(r \delta \tau-s_1)} \left(e^{\delta^d \lambda
      p_1\EXP{\|(\us+C)\cap C \|} (p_1 + p_0e^{-\delta
        \gamma(\delta)s})} - 1 \right) \  ds \ ds_1  \\
  & \sim & n\del^2 2 e^{-2\rho p_1\beta}\int_0^{r\tau} \int_0^{r\tau - s_1}
  \left(e^{\rho p_1\EXP{\|(\us+C)\cap C \|}(p_1 + p_0e^{-a_0s})} - 1
  \right)\ ds \  ds_1. \no
\end{eqnarray*}
The last relation is obtained by replacing $s_1$ by $\delta s_1.$
Since $\EXP{|U_i - \EXP{U_i}|^3 } \leq \|D_i^1\|\VAR{U_i} = r \del
\tau \VAR{U_i}$ and $n = O(\delta^{-1})$, as $\delta \to 0$,
\begin{eqnarray}
  \frac{\sum_i\EXP{|U_i - \EXP{U_i}|}^3} {\left( \sum_i\VAR{U_i}
    \right)^{3/2}} \ \leq \ \frac{r \del \tau }{\left(
      \sum_i\VAR{U_i} \right)^{1/2}} = O(\delta^{1/2}) \to 0.
  \label{Lyapunov_conditon}
\end{eqnarray}
Hence the first part of (\ref{eqn:for_unrel_clt}) follows from
Lyapunov's central limit theorem and the second part of the proof
follows as in \cite{Hall88} (pp. 157-158).  \qed
\remove{To prove (\ref{eqn:clt2}), use (\ref{eqn_scal_var2}) to obtain an
expression for $\VAR{V^{(1)}}.$ Make the change of variable from $s_1$ to
$s_1/\del$. The estimate in (\ref{Lyapunov_conditon}), with $2(1 +
\epsilon)$ instead of $3$, will be $O(\delta
\gamma(\delta)^{\epsilon})$ which converges to $0.$}

{\bf Proof of Theorem \ref{thm:unrel_sl}:} We divide the interval
$[0,T]$ into intervals of length $2\tau/c$ where the number of odd
numbered intervals is $n$ and the number of even numbered intervals is
$m.$ $n$ and $m$ are such that $n+m = \lfloor c T/2\tau \rfloor$ and
$0 \leq n-m \leq 1.$ Let $U_i^1$ be the vacancy of the $i$-th
odd-numbered interval and $U_i^2$ the vacancy of the $i$-th
even-numbered interval. Hence,
\begin{eqnarray*}
  V_T & = & \int_0^T 1[\us \notin \cC_s]ds \\
  & = & \sum_{i=0}^n \int_{4i\tau/c}^{2(2i+1)\tau/c} 1[\us \notin \cC_s]ds +
  \sum_{i=0}^m \int_{2(2i+1)\tau/c}^{4(i+1)\tau/c} 1[\us \notin \cC_s]ds +
  \int_{(n+m)2\tau/c}^T 1[\us \notin \cC_s]ds \\
  & = & \sum_{i=0}^n U_i^1 + \sum_{i=0}^m U_i^2 + \int_{(n+m)2\tau/c}^T 1[\us
  \notin \cC_s]ds.
\end{eqnarray*}
As $T \to \infty$,
\begin{displaymath}
  \frac{1}{T}\int_{(n+m)2\tau/c}^T 1[\us \notin \cC_s]ds \ \leq \
  \frac{2\tau}{cT} \to 0.
\end{displaymath}
$\{U_i^1\}$ and $\{U_i^2\}$ are sequences of i.i.d. random variables
while $U_i^1$ and $U_i^2$ might be dependent for some $i$'s. Hence, by
the strong law of large numbers,
\begin{displaymath}
  n^{-1}\sum_{i=0}^n U_i^1  \to \EXP{U_i^1} = 2\tau e^{-p_1 \lam \beta}/c
  \qquad \mbox{a.s.},
\end{displaymath}
and
\begin{displaymath}
  m^{-1}\sum_{i=0}^m U_i^2  \to \EXP{U_i^2} = 2\tau e^{-p_1 \lam \beta}/c
  \qquad \mbox{a.s.}.
\end{displaymath}
Since
\begin{displaymath}
  \frac{V_T}{T} =\left( \frac{n}{T}\right) n^{-1}\sum_{i=0}^n U_i^1  +
  \left( \frac{m}{T}\right) m^{-1} \sum_{i=0}^m U_i^2 +
  \left( \frac{1}{T} \right)\int_{(n+m)2\tau/c}^T 1[\us \notin
  \cC_s] \ ds,
\end{displaymath}
the strong law follows by noting that $n/T, m/T \to 1/4\tau$ as $T \to
\infty.$ \qed

{ \bf Proof of Theorem \ref{thm:unrel_clt2}:} The proof follows using
a similar spatial blocking technique as in Theorem \ref{thm:unrel_clt}
provided we show $T^{-1} \ \VAR{V_T} \rar \sigma_2^2$ as $T \rar
\infty.$
First make the change of the variable $s=s_2-s_1$ and then $u =
s_1/T$ in (\ref{eqn:var(V_T)}), to obtain
\[ \VAR{V_T} = 2Te^{-2 p_1 \lam \beta} \int_0^1 \int_0^{T(1-u)}
\left(e^{p_1 \lam \EXP{\|(\us+C)\cap C\|}(p_1 + p_0e^{-(\mu_0 +
      \mu_1)s})} - 1 \right) \ ds \ du. \]
The result now follows from the monotone convergence theorem and using the fact that $C \subset B_{0,\tau}$.
The finiteness of the limiting variance follows from uniform boundedness of inner integral in the above equation.  Indeed for sufficiently large $T,$ we have
\begin{eqnarray*}
\int_0^{T(1-u)}
\left(e^{p_1 \lam \EXP{\|(\us+C)\cap C\|}(p_1 + p_0e^{-(\mu_0 +
      \mu_1)s})} - 1 \right) \ ds & \leq & \int_0^{T}
\left(e^{p_1 \lam \EXP{\|(\us+C)\cap C\|}(p_1 + p_0e^{-(\mu_0 +
      \mu_1)s})} - 1 \right) \ ds \\
 & = & \int_0^{2\tau/c}
\left(e^{p_1 \lam \EXP{\|(\us+C)\cap C\|}(p_1 + p_0e^{-(\mu_0 +
      \mu_1)s})} - 1 \right) \ ds \\
 & \leq & \frac{2 \tau}{c} \left( e^{\lam \EXP{\|C\|}} + 1 \right) < \infty.
\end{eqnarray*}
\qed

\section{Appendix}
\label{sec:app}
{ \bf Proof of Inequality \ref{coverage_bounds}: }
Let $R$ be the unit square $[0,1]^2$. Consider the PB model defined in Subsection~\ref{sec:Markov-Boolean} with intensity $\lambda$ and the sensing regions to be discs of radius $r > 0,$ satisfying $\pi r^2 < 1.$ Let $V_k$ (see (\ref{eq:V_k-defn})) be the $k-$vacancy in the region $R$.  First we shall derive the upper bound. We can write
\begin{equation}
\prob{V_k > 0} = p_1 + p_2 + p_3 ,
\label{p0}
\end{equation}
where
\begin{eqnarray}
p_1 & = & \prob{\mbox{no disk is centered within $R$}} = e^{-\lam} \leq e^{-\lam \pi r^2}, \qquad \mbox{since $\pi r^2 < 1$},
\label{p1} \\
p_2 & = & \mathsf{Pr}\left(\mbox{at least one disk is centered within $R$ but none of the disks intersect any other disk and} \right. \nonumber \\
& & \left. \mbox{none of the disks intersect the boundary of $R$} \right) \nonumber \\
& \leq & \prob{\mbox{at least one disk is centered within R}} \times \prob{\mbox{a given disk intersects no other disk}}
\nonumber \\
& = & (1 - e^{-\lam}) \times e^{-\lam \pi (2r)^2} \leq e^{-\lam \pi r^2}, \qquad \mbox{and} \label{p2} \\
p_3 & = & \mathsf{Pr} \left( \mbox{$R$ is not completely $k$-covered, at least one disk is centered
within $R$ and} \right. \nonumber  \\
& & \left. \mbox{at least two disks intersect each other or at least one disk intersects the boundary of $R$}\right).
\nonumber
\end{eqnarray}
We shall now obtain an upper bound for $p_3$ using the same technique as in \cite{Zhang04}. Define a \emph{crossing}
to be either the point of intersection of two disks or the point of intersection of a disk with the boundary of the region $R$.  A crossing is said to be $k$-covered if it is an interior point of at least $k$ disks. By Theorem 4 of \cite{Wang03}, $R$ is completely $k$-covered iff there exist crossing points and every crossing point is $k$-covered.

Let $N$ and $M$ be the number of crossings and the number of crossings that are not $k$-covered respectively. Then, by definition of crossings and its relation to complete $k$-coverage, we have
\begin{eqnarray}
p_3 & \leq & \prob{M \geq 2}  \leq \frac{\EXP{M}}{2} \nonumber \\
 & \leq & \half \EXP{N}\prob{\mbox{a crossing is not $k$-covered}} \nonumber  \\
 & = & \half (\EXP{N_1}+\EXP{N_2})\left(e^{-\lam \pi r^2}\sum_{i=0}^{k-1}\frac{(\lam \pi r^2)^i}{i!} \right), \label{p3a}
\end{eqnarray}
where $N_1$ is the number of crossings created by intersection of two disks and $N_2$ is the number of crossing created by intersection of the disks with the boundary of $R$. By Palm theory for Poisson processes (see \cite{Stoyan95}), we have the following
\begin{equation}
\EXP{N_1}  =  \EXP{\mbox{No. of pts in $R$}}\EXP{\mbox{No. of crossings created by a given node with other nodes}} = 4\lam^2 \pi r^2. \label{n1}
\end{equation}
In the above formula, the two crossings created by intersection of two disks are assigned one each to the two nodes. For estimating $\EXP{N_2}$ note that only nodes within a distance of $r$ from the boundary can create crossings with the boundary. The expected number of such nodes is at most $4\lam r$ and each node can create atmost $2$ crossings. Again by Palm theory, we have
\begin{equation}
\EXP{N_2} \leq  8 \lam r. \label{n2} \end{equation}
Substituting from (\ref{n1}), (\ref{n2}) in (\ref{p3a}), we get

\[ p_3 \leq 2 \lam^2 \pi r^2(1 + \frac{2}{\lam \pi r}) e^{-\lam \pi r^2}\sum_{i=0}^{k-1}\frac{(\lam \pi r^2)^i}{i!}.\]
The upper bound in (\ref{coverage_bounds}) now follows from the above inequality and (\ref{p0})-(\ref{p2}). The lower bound is based on the following inequality which follows from (3.7) of \cite{Hall88}.
\begin{equation}
\prob{V_k > 0} \geq \frac{(\EXP{V_k})^2}{\EXP{V_k^2}}. \label{lb_V_k}
\end{equation}
By definition,
\[ \EXP{V_k^2}  =  I_1 + I_2, \]
where
\begin{eqnarray*}
I_1 & = & \int \int_{R^2 \cap \{|x_1-x_2| > 2r\}}\EXP{V_k(x_1)V_k(x_2)} \ dx_1 \ dx_2 \leq (\EXP{V_k})^2, \\
I_2 & = & \int \int_{R^2 \cap \{|x_1-x_2| \leq 2r\}}\EXP{V_k(x_1)V_k(x_2)} \ dx_1 \ dx_2.
\end{eqnarray*}
The bound for $I_1$ above is obtained using the independence of the two terms in the integrand.
Let $B_{x_i, r}$ be disks of radius $r$ centered at $x_i ,$ $i=1,2$ respectively.
By standard calculations (see proof of Theorem 3.11, \cite{Hall88}), it can be shown that $\|B_{x_2,r} \setminus B_{x_1,r}\| \geq \pi r |x_1 - x_2|/2$.
\begin{eqnarray*}
I_2 & = &  \int \int_{R^2\cap\{|x_1-x_2|\leq 2r\}}\EXP{V_k(x_1)V_k(x_2)} \ dx_1 \ dx_2 \\
& = &\int \int_{R^2\cap\{|x_1-x_2|\leq 2r\}} \prob{\mbox{both $B_{x_1,r}$ and $B_{x_2,r}$ contain less than $k$ nodes}} \ dx_1 \ dx_2 \\
& \leq & \int \int_{R^2\cap\{|x_1-x_2|\leq 2r\}} 
\prob{\mbox{\small{ $B_{x_1,r}$ contains less than $k$ nodes,$B_{x_2,r} \setminus B_{x_1,r}$ contains less than $k$ nodes}}} \ dx_1 \ dx_2 \\
& = & \int \int_{R^2\cap\{|x_1-x_2|\leq 2r\}} \prob{\mbox{\small{$B_{x_1,r}$ contains less than $k$ nodes}}}
\prob{\mbox{\small{$B_{x_2,r} \! \setminus \! B_{x_1,r}$ contains less than $k$ nodes}}}  dx_1 dx_2 \\
& \leq & \int \int_{R^2\cap\{|x_1-x_2|\leq 2r\}} \EXP{V_k(x_1)} \left(e^{-\lam \pi r|x_1-x_2|/2}\sum_{i=0}^{k-1}\frac{(\lam \pi r|x_1-x_2|/2)^i}{i!} \right)  dx_1 dx_2 \\
& = & \int_R \EXP{V_k(x_1)} \ dx_1 \int_0^{2r} e^{-\lam \pi rx/2}\sum_{i=0}^{k-1}\frac{(\lam \pi rx/2)^i}{i!} 2 \pi x \ dx \\
& = & \frac{\EXP{V_k}}{\pi r^2} \int_0^{\pi r^2} e^{-\lam u}\sum_{i=0}^{k-1}\frac{(\lam u)^i}{i!}8u \ du \\
& \leq & \frac{\EXP{V_k}}{\pi r^2} \int_0^{\infty} e^{-\lam u}\sum_{i=0}^{k-1}\frac{(\lam u)^i}{i!}8u \ du = \EXP{V_k} 4k(k+1)\frac{\lam^{-2}}{\pi r^2},\\
\end{eqnarray*}
where the equality in the penultimate line is obtained by making the change of variable $u = \pi r x/2$. The last equality follows from the usual representation of the Gamma function (see (26), \cite{Zhang04}).
Substituting the above bounds for $I_1$ and $I_2$ in (\ref{lb_V_k}), we get
\begin{equation}
\prob{V_k > 0} \geq \frac{1}{1 + \frac{4k(k+1)}{\lam^2 \pi r^2 \EXP{V_k}}}.
\label{bound_prob_v_k_pos}
\end{equation}
Setting $\beta = \pi r^2$ in (\ref{eq:expected-V_k}) and ignoring all but the last term, we get
\begin{equation}
\EXP{V_k} \geq \frac{e^{-\lam \pi r^2} (\lam \pi r^2)^{k-1}}{(k-1)!}.
\label{exp_v_k_bound}
\end{equation}
The lower bound in (\ref{coverage_bounds}) now follows from (\ref{bound_prob_v_k_pos}), (\ref{exp_v_k_bound}). \qed

{\bf Acknowledgments. } The authors would like to thank an anonymous referee for numerous suggestions leading to significant improvements in the paper.

\bibliographystyle{siam}

\bibliography{coverage}

\end{document}